\documentstyle[12pt,english]{article} 
\input amssym.def 
\input amssym.tex

\def \cst {\mbox{${\rm C}^*$}}
\def \K{{\cal K}}
\def \E {{\cal E}}
\def \D {{\cal D}}
\def \L {{\cal L}}
\def \O2{{{\cal O}_2}}
\def \Oinf{{{\cal O}_\infty}}
\def \N{{\Bbb N}}

\def \C{{\Bbb C}}
\def \R{{\Bbb R}}
\def \U{{\cal U}}
\def \ot{\otimes}

\def \mapright#1{\smash{\mathop{\longrightarrow}\limits^{#1}}}

\newtheorem{lemme}{Lemma}[section]
\newtheorem{vth}[lemme]{Theorem}
\newtheorem{defi}[lemme]{Definition}
\newtheorem{prop}[lemme]{Proposition}
\newtheorem{cor}[lemme]{Corollary}

\title{{\bf Subtriviality of continuous fields of nuclear \cst-algebras}}
\author{By Etienne Blanchard}
\date{}

\begin{document}
\maketitle
\begin{abstract} 
We extend in this paper the characterisation 
of a separable nuclear \cst-algebra given by Kirchberg proving that 
given a unital separable continuous field of nuclear \cst-algebras $A$ 
over a compact metrizable space $X$, the $C(X)$-algebra $A$ 
is isomorphic to a unital $C(X)$-subalgebra of the trivial continuous 
field $\O2\ot C(X)$, image of $\O2\ot C(X)$ by a norm one projection.

\begin{center}
\end{center}
\end{abstract}

\setcounter{section}{-1}
\section{Introduction} 
\indent\indent 
In order to study deformations in the \cst-algebraic framework, Dixmier 
introduced the notion of continuous field of \cst-algebras over a 
locally compact space (\cite{di}). 
In the same way as there is a faithful representation in a Hilbert 
space for any \cst-algebra thanks to the Gelfand--Naimark--Segal 
construction, a separable continuous field of \cst-algebras $A$ over 
a compact metrizable space $X$ always 
admits a continuous field of faithful representations $\pi$ in a Hilbert 
$C(X)$-module, i.e. there exists a family of representations $\{\pi_x$, 
$x\in X\}$, in a separable Hilbert space $H$ which factorize through a 
faithful representation of the fibre $A_x$ such that for each $a\in A$, 
the map $x\mapsto\pi_x(a)$ is strongly continuous 
(\cite[th\'eor\`eme 3.3]{bla1}). 

In a work on tensor products over $C(X)$ of continuous fields of 
\cst-algebras over $X$ (\cite{kirch}), 
Kirchberg and Wassermann raised the question of 
whether the continuous field of \cst-algebras $A$ could 
be subtrivialized, i.e. whether one could 
find a continuous field of faithful 
representations $\pi$ such that the map $x\mapsto\pi_x(a)\in L(H)$ 
is actually norm continuous for all $a$ in $A$. 
In this case, given any \cst-algebra $B$, the minimal tensor product 
$A\ot B$ is a $C(X)$-subalgebra of the trivial continuous field 
${[L(H)\ot B]\ot C(X)}$ and is therefore a continuous field with fibres 
$(A\ot B)_x=A_x\ot B$. 
They proved that a non-exact continuous field with exact 
fibres cannot be subtrivialized and they constructed such examples. 

\medskip
The non-trivial example of the continuous field of rotation algebras over 
the unit 
circle ${\Bbb T}$ had already been studied by Haagerup and R\o rdam 
in \cite{haagro}. More precisely, they constructed continuous functions 
$u,v$ from ${\Bbb T}$ to the unitary group $U(H)$ of the 
infinite-dimensional separable Hilbert space $H$ satisfying the 
commutation 
relation ${u_tv_t=tv_tu_t}$ for all $t\in {\Bbb T}$ and 
the uniform continuity condition 
$\max\left\{ \| u_t-u_{t'}\| ,\| v_t-v_{t'}\|\right\} <C'|t-t'|^{1/2}$ 
where $C'$ is a computable constant. 

Our purpose in the present paper is to show that the subtrivialization 
is always possible in the nuclear separable case through a generalisation 
of the following theorem of Kirchberg using ${\cal R}KK$-theory arguments: 
\begin{vth}\label{kirchb} (\cite{kir2}) 
A unital separable \cst-algebra $A$ is exact 
if and only if it is isomorphic to a \cst-subalgebra of $\O2$. 
Moreover the \cst-algebra $A$ is nuclear if and only if 
$A$ is isomorphic to a \cst-subalgebra of $\O2$ containing the unit 
$1_\O2$ of $\O2$, image of $\O2$ by a unital completely 
positive projection. 
\end{vth}
As a matter of fact, we get an equivalent characterisation of nuclear 
separable continuous fields of \cst-algebras (theorem \ref{concl}) 
which is made possible thanks to $C(X)$-linear homotopy invariance 
of the bifunctor ${\cal R}KK(X;--,--)$ (theorem \ref{homoto}) 
and $C(X)$-linear Weyl-von Neumann absorption results (proposition 
\ref{result}). This also enables us to have a better understanding of the 
characterisation of separable continuous fields of nuclear \cst-algebras 
given by Bauval in \cite{bauv}. 

\medskip
In an added appendix, the corresponding characterisation of 
exact separable continuous fields of \cst-algebras as $C(X)$-subalgebras 
of $\O2\ot C(X)$ given by Eberhard Kirchberg is described 
(theorem \ref{exact}). 

\medskip 
{\sl I would like to thank E. Kirchberg for his 
enlightenment on the exact case. I also want to express my gratitude to 
C. Anantharaman-Delaroche and J. Cuntz for fruitful discussions.}
\section{Preliminaries}
\subsection{$C(X)$-algebras} 
Let $X$ be a compact Hausdorff space and $C(X)$ be the \cst-algebra of 
continuous functions on $X$ with complex values. 
We start by recalling the following definition. 
\begin{defi} (\cite{kas3}) 
A $C(X)$-algebra is a \cst-algebra $A$ endowed 
with a unital morphism from $C(X)$ in the centre of the multiplier 
algebra $M(A)$ of $A$. 
\end{defi}
{\mbox{{\bf Remark:}}} We do not assume that $C(X)$ 
embeds into $M(A)$. For instance, there is a natural structure of 
$C([0,2])$-algebra on the \cst-algebra $C([0,1])$. 

\medskip 
For $x\in X$, define the kernel $C_x(X)$ of the evaluation map 
$ev_x :C(X)\rightarrow\C$ at $x$; denote by $A_x$ the quotient 
of a $C(X)$-algebra $A$ by the closed ideal $C_x(X)A$ and 
by $a_x$ the image of an element $a\in A$ in the fibre $A_x$. 
Then the function 
\begin{center}
$x\mapsto\| a_x\| =\inf\{\| [1-f+f(x)]a\| ,f\in C(X)\}$ 
\end{center}
is upper semi-continuous for any $a\in A$ 
and the $C(X)$-algebra $A$ is said to be a continuous field 
of \cst-algebras over $X$ if the function $x\mapsto\| a_x\|$ is actually 
continuous for every $a\in A$ (\cite{di}). 

\medskip\noindent {\bf Examples} 
\noindent {\sl 1.} If $A$ is a $C(X)$-algebra and $D$ is a \cst-algebra, 
the spatial tensor product $B=A\ot D$ is naturally endowed with a 
structure of $C(X)$-algebra through the map 
$f\in C(X)\mapsto f\ot 1_{M(D)}\in M(A\ot D)$. 
In particular, if $A=C(X)$, the tensor product $B$ 
is a trivial continuous field over $X$ with constant fibre $B_x\simeq D$

\noindent {\sl 2.} Given a $C(X)$-algebra $A$, define the unital 
$C(X)$-algebra ${\cal A}$ generated by $A$ and $u[C(X)]$ in 
$M[A\oplus C(X)]$ where $u(g)(a\oplus f)=ga\oplus gf$ for $a\in A$ 
and $f,g\in C(X)$. It defines a continuous field of \cst-algebras 
over $X$ if and only if the $C(X)$-algebra $A$ is continuous 
(\cite[proposition 3.2]{bla1}). 

\medskip\noindent 
{\mbox{{\bf Remark:}}} 
If $A$ is a separable continuous field of non-zero 
\cst-algebras (not necessarily unital) over the compact Hausdorff space 
$X$, the positive cone $C(X)_+$ and so the \cst-algebra $C(X)$ are 
separable. Hence, the topological space $X$ is metrizable. 

\medskip
\begin{defi} (\cite{bla1,bla2}) 
Given a continuous field of \cst-algebras $A$ over the compact Hausdorff 
space $X$, 
a continuous field of representations of a $C(X)$-algebra $D$ 
in the multiplier algebra $M(A)$ of $A$ is a $C(X)$-linear morphism 
$\pi :D\rightarrow M(A)$, i.e. for each $x\in X$, 
the induced representation $\pi_x$ of $D$ in $M(A_x)$ factorizes through 
the fibre $D_x$. 
\end{defi}

If the $C(X)$-algebra $D$ admits a continuous field of faithful 
representations $\pi$ in the $C(X)$-algebra $M(A)$ where $A$ is 
a continuous field over $X$, i.e. the induced representation 
of the fibre $D_x$ in $M(A_x)$ is faithful for every point $x\in X$, 
the function 
\begin{center}
$x\mapsto \|\pi_x(d)\| =\sup\{\| (\pi (d)a)_x\| , 
a\in A \,\mbox{such that}\, \| a\|\leq 1\}$ 
\end{center}
is lower semi-continuous for all $d\in D$ and 
the $C(X)$-algebra $D$ is therefore continuous.

In particular a separable $C(X)$-algebra $D$ is 
continuous if and only if there exists a Hilbert $C(X)$-module $\E$ 
such that $D$ admits a continuous field of faithful representations 
in the multiplier algebra $M(\K (\E) )=\L (\E )$ of the continuous field 
over $X$ of compact operators $\K (\E )$ acting on $\E$ 
(\cite[th\'eor\`eme 3.3]{bla1}). 

\bigskip\bigskip 
Let us also mention the characterisation of separable continuous fields 
of nuclear \cst-algebras over a compact metrizable space $X$ 
given by Bauval in \cite{bauv} using a natural 
$C(X)$-linear version of nuclearity introduced by Kasparov and Skandalis 
in \cite{kask}$\S$6.2 : 
a $C(X)$-linear completely positive $\sigma$ from a $C(X)$-algebra $A$ 
into a $C(X)$-algebra $B$ is said to be $C(X)$-nuclear if and only if 
given any compact set 
$F$ in $A$ and any strictly positive real number $\varepsilon$, 
there exist an integer $k$ and $C(X)$-linear completely 
positive contractions $T: A\rightarrow M_k(\C )\ot C(X)$ and 
$S: M_k(\C )\ot C(X)\rightarrow B$ such that for all $a\in F$, 
one has the inequality 
\begin{center} 
$\| \sigma (a)-(S\circ T)(a)\| <\varepsilon$. 
\end{center}
One can then state the following results. The first assertion is a simple 
$C(X)$-linear reformulation of the Choi-Effros theorem and the second 
one is due to Bauval. 
\begin{prop}\label{chnuc} Let $X$ be a compact metrizable space and 
$A$ be a separable $C(X)$-algebra. 
\begin{enumerate}
\item (\cite{kask}$\S$6.2) 
Given a $C(X)$-algebra $B$ and a closed ideal $J\subset B$, 
any contractive $C(X)$-nuclear map $A\rightarrow B/J$ admits a 
contractive $C(X)$-linear completely positive lift $A\rightarrow B$. 
\item (\cite[th\'eor\`eme 7.2]{bauv}) 
The $C(X)$-algebra $A$ is 
a continuous fields of nuclear \cst-algebras over $X$ if and only if the 
identity map $id_A: A\rightarrow A$ is $C(X)$-nuclear. 
\end{enumerate} 
\end{prop}
{\bf Remark:} In assertion {\sl 1.}, the ideal $J=(C(X)B)J=C(X)J$ is 
a $C(X)$-algebra. 
\subsection{$C(X)$-extensions} 
Given a compact Hausdorff space $X$, we introduce a natural 
$C(X)$-linear version of the semi-group $Ext(-,-)$ defined by Kasparov 
(\cite{kas2,kas3}). 

\medskip
Call a morphism of $C(X)$-algebras a $*$-homomorphism between 
$C(X)$-algebras which is $C(X)$-linear. 
\begin{defi} 
A $C(X)$-extension of a 
$C(X)$-algebra $A$ by a $C(X)$-algebra $B$ is a short exact sequence 
\begin{center}
$0\rightarrow B\rightarrow D\mapright{\pi} A\rightarrow 0$ 
\end{center}
in the category of $C(X)$-algebras. The $C(X)$-extension is said to be 
trivial if the map $\pi$ admits a cross section 
$s:A\rightarrow D$ which is a morphism of $C(X)$-algebras. 
\end{defi}

As in the \cst-algebraic case a $C(X)$-extension 
$0\rightarrow B\rightarrow D\rightarrow A\rightarrow 0$ of 
$A$ by $B$ defines unambiguously an homomorphism from $D$ 
to the multiplier algebra $M(B)$ of $B$, which gives 
a morphism of $C(X)$-algebras $\sigma :A\rightarrow M(B)/B$ 
(called the Busby invariant of the extension) 
and the $C(X)$-extension is trivial if and only if the map $\sigma$ lifts 
to a morphism of $C(X)$-algebras $A\rightarrow M(B)$. 
Conversely, given a morphism of $C(X)$-algebras 
$\sigma :A\rightarrow M(B)/B$, setting 
$D=\{ (a,m)\in A\times M(B), \sigma (a)=q(m)\}$ where $q$ is the quotient 
map $M(B)\rightarrow M(B)/B$, one has a $C(X)$-extension 
$0\rightarrow B\rightarrow D\rightarrow A\rightarrow 0$ 
(see \cite{kas2}$\S$7). 

\medskip\noindent 
{\mbox{{\bf Remark:}}} A $C(X)$-extension 
$0\rightarrow B\rightarrow D\rightarrow A\rightarrow 0$ 
induces for every $x\in X$ a \cst-extension 
$0\rightarrow B_x\rightarrow D_x\rightarrow A_x\rightarrow 0$.

\bigskip\bigskip 
In order to define the sum of two $C(X)$-extensions, recall that 
the Cuntz algebra $\O2$ is the unital \cst-algebra 
generated by two orthogonal isometries $s_1$ and $s_2$ subject 
to the relation $1=s_1s_1^*+s_2s_2^*$ (\cite{cuntz}). 
Then if $\K$ is the \cst-algebra of compact operators 
on the infinite-dimensional separable 
Hilbert space, one defines the sum of two $C(X)$-extensions 
$\sigma_1$ and $\sigma_2$ of the $C(X)$-algebra $A$ by 
the stable $C(X)$-algebra $\K\ot B$ through the choice of 
a unital copy of $\O2$ in the multiplier algebra $M(\K )$ of $\K$ 
to be the $C(X)$-extension 
\begin{center}
$\sigma_1\oplus\sigma_2 :a\mapsto
q(s_1\ot 1)\sigma_1(a)q(s_1^*\ot 1)+q(s_2\ot 1)\sigma_2(a)q(s_2^*\ot 1)
\in M(\K\ot B)/(\K\ot B)$, 
\end{center}
where $q$ is the quotient map 
$M(\K\ot B)\rightarrow M(\K\ot B)/(\K\ot B)$. 
\begin{defi} Given a compact Hausdorff space $X$ and two $C(X)$-algebras 
$A$ and $B$, $Ext(X;A,B)$ is the semi-group of $C(X)$-extensions of $A$ 
by $\K\ot B$ divided by the equivalence relation $\sim$ where 
$\sigma_1\sim\sigma_2$ if there exist a unitary $U\in M(\K\ot B)$ of 
image $q(U)$ in the quotient $M(\K\ot B)/(\K\ot B)$ and 
two trivial $C(X)$-extensions ${\pi_1}$ and ${\pi_2}$ such that 
for all $a\in A$, 
\begin{center}
$(\sigma_2\oplus {\pi_2} )(a)=
q(U)^*(\sigma_1\oplus {\pi_1} )(a)q(U)$ (in $M(\K\ot B)/(\K\ot B)$). 
\end{center}
\end{defi} 

Let $Ext(X;A,B)^{-1}$ be the group of invertible elements of 
$Ext(X;A,B)$, i.e. the group of classes of $C(X)$-extension 
$\sigma$ such that there exists a $C(X)$-extension $\theta$ with 
$\sigma\oplus\theta$ trivial. 
One can generalise Kasparov's theorem of homotopy invariance 
of the group $Ext (A,B)^{-1}$ to the framework of $C(X)$-algebras as 
follows. 
\begin{vth}\label{homoto}(\cite{kas2}) 
Assume that $A$ is a separable $C(X)$-algebra and that 
$B$ is a $\sigma$-unital $C(X)$-algebra. 
Then the group $Ext(X;A,B)^{-1}$ is isomorphic to the group 
${\cal R}{KK}^{1}(X;A,B)$ and is therefore 
$C(X)$-linear homotopy invariant in both entries $A$ and $B$. 
\end{vth}{\mbox{{\bf Proof :}}} 
Let us first make the following observation. 
Given a $C(X)$-algebra $B$ and a Hilbert $B$-module $\E$, denote 
by $\L (\E )$ the set of bounded 
$B$-linear operators on ${\E}$ which admit an adjoint (\cite{kas1}). 
Then any operator $T\in\L (\E )$ is $B$-linear and so $C(X)$-linear. 
This argument provides a natural extension of the Stinespring-Kasparov 
theorem (\cite{kas2}) to the framework of $C(X)$-algebras. 
Consequently, if $A$ is a separable 
$C(X)$-algebra and $B$ is a $\sigma$-unital $C(X)$-algebra, the class 
of a $C(X)$-extension $\sigma :A\rightarrow M(\K\ot B)/(\K\ot B)$ 
is invertible in $Ext(X;A,B)$ if and only if there is a $C(X)$-linear 
completely positive 
contractive lift $A\rightarrow M(\K\ot B)$.

\medskip
Let ${\cal R}{\Bbb E}(X;A,B)$ be the set of Kasparov 
$C(X)$--$A,B$-bimodules (\cite{kas3}definition 2.19), i.e. the set 
of Kasparov $A,B$ bimodules $(\E ,F)$ such that the representation 
$A\rightarrow\L (\E )$ is a $C(X)$-representation. 
Call a $C(X)$-linear operator 
homotopy an element $\{ (\E ,F_t),{ 0\leq t\leq 1\} }\in 
{\cal R}{\Bbb E}(X;A,B\ot C([0,1]))$ such that $t\mapsto F_t$ is norm 
continuous and define on ${\cal R}{\Bbb E}(X;A,B)$ the equivalence 
relation corresponding to the one defined by Skandalis 
in \cite[definition 2]{skand1}. 
The constructions given by Kasparov in \cite[section 7]{kas2} imply 
that, if the $C(X)$-algebra $B$ is $\sigma$-unital, 
the group of equivalence classes 
${\cal R}\widetilde{KK}(X;A,B\ot {\cal C}_1)$ is 
isomorphic to $Ext(X;A,B)^{-1}$, where ${\cal C}_1$ is the first (graded) 
Clifford algebra. 

\medskip 
On the other hand, given two graded $C(X)$-algebras $A$ and $B$ 
with $A$ separable, the different steps of the demonstration of 
\cite[theorem 19]{skand1} provide us with an isomorphism between 
the two groups ${\cal R}\widetilde{KK}(X;A,B)$ and ${\cal R}{KK}(X;A,B)$ 
since proposition 2.21 of \cite{kas3} defines 
an intersection product in ${\cal R}\widetilde{KK}$-theory 
and lemma 18 of \cite{skand1} gives us the equality 
\begin{center}
$(ev_0\ot id_{C(X)})^*(1_{C(X)} )=(ev_1\ot id_{C(X)})^*(1_{C(X)} )$ 
in ${\cal R}{\widetilde KK}(X;C([0,1])\ot C(X),C(X))$, 
\end{center}
where $1_{C(X)}$ is the Kasparov $C(X),C(X)$-bimodule $(C(X),0)$ and 
$ev_t:C([0,1])\rightarrow\C$ is the evaluation map at $t\in [0,1]$. 
$\square$

\medskip\noindent {\bf Remarks:} {\sl 1.} Kuiper's theorem implies 
that the law of addition on the abelian group $Ext(X;A,B){}^{-1}$ 
is independent of the choice of the unital copy of $\O2$ in $M(\K )$. 

\medskip\noindent 
{\sl 2.} If $A$ is a separable nuclear continuous field of \cst-algebras 
over $X$ and $B$ is a $C(X)$-algebra, every $C(X)$-linear morphism 
from $A$ to the quotient $M(\K\ot B)/(\K\ot B)$ is $C(X)$-nuclear 
and therefore admits a $C(X)$-linear completely 
positive lifting $A\rightarrow M(\K\ot B)$ thanks to 
proposition \ref{chnuc}. Accordingly one has the equality 
\begin{center}
$Ext(X;A,B){}^{-1}=Ext(X;A,B )$. 
\end{center}

\section{An absorption result}
In this section we prove a continuous generalisation 
of a statement contained in \cite{kir2} which will enable us 
to get a $C(X)$-linear Weyl-von Neumann type result 
(proposition~\ref{result}). 
Let us start with the following definition of Cuntz. 
\begin{defi} (\cite{cuntz}) 
A simple \cst-algebra distinct from $\C$ is said to be purely infinite 
if and only if for any non-zero $a,b\in A$, there exist elements 
$x,y\in A$ such that $a=xby$. 
\end{defi} 
Then, we can state a proposition from Kirchberg's classification work, 
based on Glimm's lemma (\cite{di}, $\S$ 11.2). A sketch of proof can also 
be found in \cite[proposition 5.1]{claire}. 
\begin{prop}(\cite{kir2}) Let $A$ be a purely infinite simple 
\cst-algebra and 
$D$ be a separable \cst-subalgebra of $M(A)$. 
Assume that $V:D\rightarrow A$ is a nuclear contraction. 

Then there exists a sequence $(a_n)$ of elements in $A$ of norm less 
than $1$ such that $V(d)=\lim_{n\rightarrow\infty} a_n^*da_n$ for 
all $d\in D$. 
\end{prop} 
{\bf Remark:} A simple ring has by definition exactly 
two distinct two sided ideals and is therefore non-zero. 

\begin{cor}\label{techn} 
Let $A$ be a continuous field of purely infinite 
simple \cst-algebras over a compact Hausdorff space $X$ 
and assume that $D$ is a separable 
$C(X)$-subalgebra of the multiplier algebra $M(A)$ such that 
there is a unital $C(X)$-embedding of the $C(X)$-algebra 
${\cal O}_{\infty}\otimes C(X)$ in the commutant 
$D'$ of $D$ in $M(A)$ and the identity map 
$id_D: D\rightarrow M(A)$ is a continuous field of 
faithful representations. 

If $V:D\rightarrow A$ is a $C(X)$-nuclear contraction, there exists a 
sequence $(a_n)$ in the unit ball 
of $A$ with the property that for all $d\in D$, 
$$V(d)=\lim_{n\rightarrow\infty} a_n^*da_n.$$
\end{cor}{\mbox{{\bf Proof :}}} 
If $F$ is a compact generating set for $D$, it is enough to prove that 
given a strictly positive real number $\varepsilon >0$, 
there exists an element $a$ in the unit ball of $A$ such that 
$\| V(d)-a^*da\| <\varepsilon$ for all $d\in F$. 

\medskip
For $x\in X$, the fibre $A_x$ is a purely infinite simple \cst-algebra 
and the map $d\mapsto V(d)_x\in A_x$ factorizes 
through $D_x\simeq (id_D)_x(D)\subset M( A_x)$ since $id_D$ is a 
continuous field of faithful representations. As a consequence, 
the previous proposition implies that we can find an 
element $g\in A$ with $\| g\|\leq 1$ satisfying for all $d\in F$ 
the inequality 
\begin{center} 
$\| \Bigl[ V(d)-g^*dg\Bigr]_x\| <\varepsilon$. 
\end{center}
Thus, by upper semi-continuity and compactness, there exist a finite 
open covering $\{ U_1,\ldots ,U_n\}$ of the space $X$ 
and elements $g_1,\ldots ,g_n$ in the unit ball of $ A$ such that 
for all $d\in F$ and $x\in U_i$, $1\leq i\leq n$, 
\begin{center}
$\| \Bigl[ V(d)-g_i^*dg_i\Bigr]_x\| <\varepsilon$. 
\end{center}

Choose $n$ orthogonal isometries $w_1,\ldots ,w_n$ in the \cst-algebra 
$\Oinf\ot 1_{C(X)}\subset D'$ 
and let $\{\phi_i\}$ be a partition of the unit $1_{C(X)}$ subordinate 
to the covering $\{ U_i\}$ of $X$. 
The element $a=\sum_i\phi_i^{1/2}w_ig_i\in A$ verifies: 
\begin{enumerate}
\item $a^*a=\sum_{i,j}\sqrt{\phi_i\phi_j}\, g_i^*w_i^*w_jg_j 
=\sum_i\phi_ig_i^*g_i\leq 1_{M(A)}$,
\item for $d\in F$ and $x\in X$, 
$\|\Bigl[ V(d)-a^*da\Bigr]_x\|
\leq \sum_i\phi_i (x)\|\Bigl[ V(d)-g_i^*dg_i\Bigr]_x\| <\varepsilon$. 
$\square$ 
\end{enumerate}

\medskip
Let us mention the following technical corollary which will be needed in 
theorem~\ref{concl}. 
\begin{cor}\label{unite} 
If $p\in\O2\ot C(X)$ is a 
projection such that for all points $x\in X$, $p_x$ is non-zero, 
then there exists an isometry $u\in\O2\ot C(X)$ such that $p=uu^*$. 
\end{cor}{\mbox{{\bf Proof :}}} 
Let ${\cal D}_2=\lim_{n\rightarrow\infty} \O2^{\ot n}$ be the infinite 
tensor product of $\O2$. 

Given a projection $q\in {\cal D}_2\ot C(X)$ such that $\| q_x\| =1$ 
for all $x\in X$, we first show that there exists an element 
$v\in {\cal D}_2\ot C(X)$ satisfying $1_{{\cal D}_2\ot C(X)}=v^*qv$. 
Namely, by density of the algebraic tensor product 
\begin{center}
$\left[ \mathop{\cup}\limits_n \O2^{\ot n}\right]\odot C(X)=
\mathop{\cup}\limits_n \left[\O2^{\ot n}\odot C(X)\right]$ 
\end{center}
in the \cst-algebra ${\cal D}_2\ot C(X)$ and functional calculus 
one can find an integer $n>0$ and a projection 
$r\in\O2^{\ot n}\ot C(X)\subset {\cal D}_2\ot C(X)$ 
such that $\| q-r\| <1$, which implies that $r=s^*qs$ for some element 
$s\in {\cal D}_2\ot C(X)$. 
Take then a faithful state $\varphi$ on $\O2^{\ot n}$ and consider the 
$C(X)$-linear completely positive map 
\begin{center}
$V: [\O2^{\ot n}\ot 1_\O2]\ot C(X) 
\rightarrow \O2^{\ot n+1}\ot C(X)$ 
\end{center}
defined by the formula 
$V(d)=(\varphi\ot id_{C(X)}) (d)1_{\O2^{\ot n+1}\ot C(X)}$ 
for $d\in [\O2^{\ot n}\ot 1_\O2]\ot C(X)\simeq \O2^{\ot n}\ot C(X)$. 
According to corollary \ref{techn}, there exists an element 
$t\in \O2^{\ot n+1}\ot C(X)$ such that 
$$1_{{\cal D}_2\ot C(X)}=1_{\O2^{\ot n+1}\ot C(X)}=
t^*rt=(st)^*q(st).$$

\medskip 
Consider now the set $\cal P$ of projections $p$ in $\O2\ot C(X)$ 
such that $p_x\not =0$ for all points $x\in X$. 
If $p$ belongs to $\cal P$, 
there exists an isometry $v\in \O2\ot C(X)$ such that $p\geq vv^*$ since 
the $K$-trivial purely infinite separable unital nuclear \cst-algebra 
${\cal D}_2$ satisfying the U.C.T. is isomorphic to $\O2$ (\cite{kir2}). 
As a consequence, if $t$ is the isometry $t=v(s_1\ot 1)v^*$, 
the projection $r=tt^*$ (Murray-von Neumann 
equivalent to $1_{\O2\ot C(X)}$) verifies 
\begin{center} 
$p-r\geq r'=v(s_2s_2^*\ot 1)v^*\in\cal P$. 
\end{center} 

The non-empty set $\cal P$ therefore satisfies 
the conditions \mbox{$(\pi_1)$--$(\pi_4)$} defined by Cuntz 
in \cite{cuntz}. But the \cst-algebra $\O2\ot C(X)$ is 
$K_0$-triviality thanks to \cite[theorem 2.3]{cuntz} and 
the theorem 1.4 of \cite{cuntz} enables us to conclude. $\square$ 

\bigskip 
One now deduces from corollary \ref{techn} the following 
absorption results (\cite{voi,kas2,kir2}): 
\begin{prop}\label{result} 
Let $A$ be a $\sigma$-unital continuous field of purely infinite 
simple nuclear \cst-algebras over a compact Hausdorff space $X$ 
and let $\K$ be the \cst-algebra of compact operators on the 
separable Hilbert space $H$. 
Denote by $q$ the quotient map 
$M(\K\ot A)\rightarrow M(\K\ot A)/(\K\ot A)$.

\begin{enumerate} 
\item Assume that $D$ is a unital separable $C(X)$-subalgebra of 
the multiplier algebra $M(\K\ot A)$ with same unit such that 
there is a unital $C(X)$-embedding of the $C(X)$-algebra 
${\cal O}_{\infty}\otimes C(X)$ in the commutant 
of $D$ in $M(\K\ot A)$ and the identity map $id_D$ is a continuous 
field of faithful representations of $D$ in $M(\K\ot A)$. 

\begin{enumerate}
\item If $V$ is a unital $C(X)$-linear completely positive map 
from $D$ in $M(\K\ot A)$ which is zero on the intersection 
$D\cap (\K\ot A)$, there exists a sequence of isometries $s_n$ in 
$M(\K\ot A)$ such that for every $d\in D$, 
\begin{center}
$V(d)-s_n^*ds_n\in\K\ot A$ and $V(d) =\lim_n s_n^*ds_n$. 
\end{center}
\item If $\pi$ is a unital morphism of $C(X)$-algebras from $D$ into 
$M(\K\ot A)$ which is zero on the intersection 
$D\cap (\K\ot A)$, there exists a sequence of unitaries $u_n$ in 
$M(\K\ot A)$ such that for every $d\in D$, 
\begin{center}
$(d\oplus\pi (d) )-u_n^*du_n\in\K\ot A$ and 
$(d\oplus\pi (d) )=\lim_n u_n^*du_n$. 
\end{center}
\item Let $B$ be a $C(X)$-algebra and assume that the quotient 
$D/(D\cap (\K\ot A) )$ is 
isomorphic to the $C(X)$-algebra ${\cal B}$, where ${\cal B}$ is 
the unital $C(X)$-algebra generated by $C(X)$ and $B$ in $M[B\oplus C(X)]$ 
(\cite[d\'efinition 2.7]{bla1}). 

Then, if $\pi :B\rightarrow M(\K\ot A)$ is a $C(X)$-linear 
homomorphism, there exists a unitary $U\in M(\K\ot A)$ such that 
for all $b\in B\subset M(\K\ot A)/(\K\ot A)$, 
\begin{center} 
$b\oplus (q\circ\pi )(b)=q(U)^*\, b\, q(U)$. 
\end{center} 
\end{enumerate}
\item 
Assume that the continuous field $A$ is separable and 
let $D$ be a separable $C(X)$-sub\-algebra of $M(A)$ containing $A$ 
such that the identity 
representation $D\rightarrow M(A)$ is a continuous field of faithful 
representations and there is a unital $C(X)$-embedding of the 
$C(X)$-algebra ${\cal O}_{\infty}\otimes C(X)$ in the commutant of $D$ 
in $M(A)$. Define the quotient $C(X)$-algebra $B=D/A$. 

If $\pi :\K\ot B\rightarrow M(\K\ot A)$ is a morphism of $C(X)$-algebras, 
there exists a unitary $U\in M(\K\ot A)$ such that 
for all $b\in\K\ot B\subset M(\K\ot A)/(\K\ot A)$, 
\begin{center} 
$b\oplus (q\circ\pi )(b)=q(U)^*\, b\, q(U)$. 
\end{center} 
\end{enumerate}
\end{prop}{\mbox{{\bf Proof :}}} 
{\it 1.} It derives from corollary \ref{techn} by the same method as
the one developed by Kasparov in \cite[theorem 5 and 6]{kas1}. 
Nevertheless, for the convenience of the reader we describe the different 
steps of the demonstration. 

\noindent {\it 1.a)} Let $F$ be a compact generating set for $D$ 
containing the unit $1_{M(\K\ot A)}$. Then given a real number 
$\varepsilon >0$, it is enough to find an element $a\in M(\K\ot A)$ 
such that $V(d)-a^*da\in\K\ot A$ and $\| V(d)-a^*da\| <3\varepsilon$ 
for all $d\in F$. 

Let $\{ e_n\}$ be an increasing, positive, quasicentral, 
countable approximate unit in the ideal $\K\ot A$ of the \cst-algebra 
generated by $\K\ot A+V(D)$. 
If we set $f_0=(e_0)^{1/2}$ and $f_k=(e_k-e_{k-1})^{1/2}$ for $k\geq 1$, 
we can then assume, passing to a subsequence of $(e_n)$ if necessary, 
that $\| V(d)f_k-f_kV(d)\| <2^{-k}\varepsilon$ for all $k\in\N$ 
and $d\in F$. This implies that the series $\sum_k [V(d)f_k-f_kV(d)]f_k$ 
is convergent in $\K\ot A$ and its norm is less than $\varepsilon$. 
Furthermore, the series $\sum_k[f_kV(d)f_k]$ is strictly convergent in 
$M(\K\ot A)$ for all $d\in F$ 
since $\sum_kf_k^2$ is strictly convergent to $1$. 

\medskip 
Notice now that the maps $V_k(d)=f_kV(d)f_k$ are all $C(X)$-nuclear 
since the separable continuous field $\K\ot A$ is nuclear. 
The corollary \ref{techn} therefore enables us to choose 
by induction $a_k\in\K\ot A$ satisfying the following conditions: 
\begin{enumerate}
\item $\forall d\in F\;\;\;\| V_k(d)-a_k^*da_k\| <2^{-k}\varepsilon$, 
\item $\forall d\in F,\forall l<k \;\;\;\|a_l^*da_k\| <2^{-l-k}\varepsilon$, 
\item $\sum_k a_k$ is strictly convergent toward an element $a\in M(A)$. 
\end{enumerate}
One then checks as in \cite[theorem 5]{kas1} that the limit $a$ 
satisfies the desired properties. 

\medskip\noindent {\it 1.b)} Take a compact generating $F$ for $D$ 
containing $1_{M(\K\ot A)}$ and consider the homomorphism 
$\pi '=1\ot\pi :D\rightarrow M(\K\ot (\K\ot A))\simeq M(\K\ot A)$. 
Given $\delta >0$, one can find, thank to the previous assertion, 
an isometry $s\in M(\K\ot A)$ such that 
\begin{center}
$s^*ds-\pi '(d)\in\K\ot A$ and $\|s^*ds-\pi '(d)\| <\delta$ 
for all $d\in K^*K$. 
\end{center}
As a consequence, if we fix $\varepsilon >0$, 
the choice of $\delta$ small enough gives us the inequality 
$\|pd-dp\| <\varepsilon$, and so 
$\| d-[pdp+p^\perp dp^\perp ]\| <2\varepsilon$ for all $d\in F$, 
where $p=ss^*$ and $p^\perp=1-p$. 

Define the unital map 
$\Theta :D\rightarrow M(p^\perp (\K\ot A)p^\perp )$ by the formula 
$\Theta (d)=p^\perp dp^\perp$. 
According to the stabilisation theorem of Kasparov 
(\cite[theorem 2]{kas1}), one can construct 
a unitary $w\in M(\K\ot A)$ verifying for all $d\in F$ the inequality 
\begin{center}
$\| d-w^*[\pi '(d)\oplus\Theta (d)]w\| <3\varepsilon$. 
\end{center}
To finish the demonstration, notice that the two homomorphisms 
$\pi '$ and $\pi '\oplus\pi$ are unitarily equivalent. 

\medskip\noindent {\it 1.c)} 
Consider the unital extension $\widetilde{\pi}$ of $\pi$ to ${\cal B}$. 
Then, the morphism $\widetilde{\pi}\circ q:D\rightarrow M(\K\ot A)$ 
reduces the demonstration to the previous assertion. 

\bigskip\noindent 
{\it 2.} The identity representation of the unital $C(X)$-algebra 
$\D =(\K\ot D)+C(X)\subset M(\K\ot A)$ is clearly a continuous field 
of faithful representations 
since the unital $C(X)$-representation $C(X)\rightarrow M(A)$ is 
a continuous field of faithful representations. 
Extend the map 
$\pi :\K\ot B=(\K\ot D)/(\K\ot A)\rightarrow M(\K\ot A)$ to a 
unital morphism of $C(X)$-algebras 
$\widetilde{\pi} :\D /(\K\ot A)\rightarrow M(\K\ot A)$. 
Applying assertion~{\it 1.b)} to the unital homomorphism 
$d\mapsto (\widetilde{\pi}\circ q )(d)$ 
from the $C(X)$-subalgebra $\D\subset M(\K\ot A)$ to 
the multiplier algebra $M(\K\ot A)$ 
now leads to the desired conclusion. $\square$

\section{The subtriviality}\label{subt} 
Given a separable continuous field of nuclear \cst-algebras $A$ over $X$, 
the strategy to prove the subtriviality of the $C(X)$-algebra $A$ will be 
the same as the one developed by Kirchberg in \cite{kir2} 
to prove theorem \ref{kirchb} whose main ideas of demonstration are 
also explained in \cite[Th\'eor\`eme 6.1]{claire}. 
We associate to $A$ a $C(X)$-extension by an hereditary 
\cst-subalgebra of the trivial continuous field ${\O2\ot C(X)}$ 
(proposition \ref{extens}) and then prove that after stabilisation, 
this $C(X)$-extension splits by ${\cal R}KK$-theory arguments 
(theorem \ref{concl}). 

\bigskip\noindent {\bf \ref{subt}.1} Let us construct the 
fundamental $C(X)$-extension associated to an exact separable 
continuous field of \cst-algebras. 
\noindent\begin{prop}\label{extens} 
Given a compact Hausdorff space $X$ and a non-zero separable unital 
exact $C(X)$-algebra $A$, 
there exist a unital $C(X)$-subalgebra $F$ of $\O2\ot C(X)$ with same 
unit and an hereditary subalgebra $I$ of $\O2\ot C(X)$ such that 
$I$ is an ideal in $F$ and 
the $C(X)$-algebra $A$ is isomorphic to 
the quotient $C(X)$-algebra $F/I$. 

Furthermore, if the topological space $X$ is perfect (i.e. without any 
isolated point) and the 
$C(X)$-algebra $A$ is continuous, the canonical map 
$F\rightarrow M(I)$ is a continuous field of faithful representations. 
\end{prop}{\mbox{{\bf Proof :}}} 
Thanks to the characterisation of separable exact \cst-algebras 
obtained by Kirchberg (theorem \ref{kirchb}), one may assume that the 
\cst-algebra $A$ is a \cst-subalgebra of $\O2$ containing the unit 
of $\O2$. 

Let $G\subset\O2\ot C(X)$ be the trivial continuous field $A\ot C(X)$ 
over $X$. Then the kernel of the $C(X)$-linear morphism 
$\pi :G\rightarrow A$ defined by $\pi (a\ot f)=fa$ is the ideal 
$J=C_\Delta (X\times X)G$ where $C_\Delta (X\times X)$ is the ideal 
in $C(X\times X)$ of functions which are zero on the diagonal. 
Indeed suppose that $T\in G$ verifies $\pi (T)=0$. 
Then given $\varepsilon >0$, take a finite number of elements $a_i\in A$, 
$f_i\in C(X)$ such that $\| T-\sum_i a_i\ot f_i\| <\varepsilon$; one has 
$\| T-\sum_i(1\ot f_i-f_i\ot 1)(a_i\ot 1)\| <\varepsilon + 
\|\pi (\sum_i a_i\ot f_i)\| <2\varepsilon$. 

Define then the hereditary subalgebra $I={J[\O2\ot C(X)]J}$ in 
$\O2\ot C(X)$ generated by $J$. It is a $C(X)$-algebra since 
it is closed by Cohen theorem (see e.g. \cite[proposition 1.8]{bla1}) 
and the product $(1_\O2 \ot f)(bc)=b(1_\O2 \ot f)c$ belongs to $I$ for all 
$f\in C(X)$ and $b,c\in I$. 
If we set $F=I+G$, the intersection $G\cap I$ is reduced 
by construction to the subalgebra $J$, and so we have a $C(X)$-extension 
\begin{center}
$0\rightarrow I\rightarrow F\rightarrow A\rightarrow 0$. 
\end{center}

\medskip
Assume now that the space $X$ is perfect and that the $C(X)$-algebra 
$A$ is continuous. We need to prove that the map 
$F_x\rightarrow M(I_x)$ is injective for each $x\in X$. Let 
$a\in G$ and $b\in I$ be two elements such that the sum $d=a+b\in F$ 
verifies for a given point $x\in X$ the equality 
\begin{center}
$d_xI_x+I_xd_x=0$ (in $[\O2\ot C(X)]_x\simeq\O2\ot\C$). 
\end{center}

To end the proof, we have to show that $d_x$ is zero. 
For every $f\in C_\Delta (X\times X)$,
one has $(bf)_x=-(af)_x\in J_x$, whence $b_x\in J_x$ and so $d_x\in G_x$. 
But the representation of $G_x\simeq A$ in 
$M(J_x)\simeq M(C_x(X)A)$ is injective since 
$X$ is perfect and $A$ is continuous, 
from which we deduce that $d_x=0$. $\square$ 

\medskip\noindent 
{\mbox{{\bf Remark:}}} With the previous notations, 
if the $C(X)$-algebra $A$ is nuclear and $\psi$ is a unital completely 
positive projection from $\O2$ onto $A$, the map 
$\pi\circ (\psi\ot id_{C(X)})$ is a unital $C(X)$-linear 
completely positive map from $\O2\ot C(X)$ onto the $C(X)$-subalgebra 
$A$ which is zero on the nuclear hereditary $C(X)$-subalgebra $I$. 

\bigskip\bigskip\noindent {\bf \ref{subt}.2} 
We can now state the main theorem: 
\begin{vth}\label{concl} 
Let $X$ be a compact metrizable space and $A$ be a unital separable 
$C(X)$-algebra 
with a unital embedding of the $C(X)$-algebra $C(X)$ in $A$. 

The following assertions are equivalent: 
\begin{enumerate}
\item $A$ is a continuous field of nuclear \cst-algebras over $X$; 
\item there exist a unital monomorphism of $C(X)$-algebras 
$\alpha :A\hookrightarrow\O2\ot C(X)$ and a unital $C(X)$-linear 
completely positive map $E:\O2\ot C(X)\rightarrow A$ such that 
$E\circ\alpha =id_A$. 
\end{enumerate}
\end{vth}{\mbox{{\bf Proof :}}} 
{\em 2}$\Rightarrow${\em 1} By assumption the identity map 
$id_A =E\circ id_{\O2\ot C(X)}\circ\alpha :A\rightarrow A$ is 
nuclear since the \cst-algebra $\O2\ot C(X)$ is nuclear and so 
the \cst-algebra $A$ is nuclear. 
Besides the $C(X)$-algebra $A$ is isomorphic to the $C(X)$-subalgebra 
$\alpha (A)$ of the continuous field $\O2\ot C(X)$ and is therefore 
continuous. 

\bigskip\noindent 
{\em 1}$\Rightarrow${\em 2} $\bullet$ Let us first deal with the case 
where the space $X$ is perfect. 

Given a unital nuclear separable continuous fields $A$ over $X$ which is 
unitaly embedded in the \cst-algebra $\O2$, consider the $C(X)$-extension 
\begin{center}
$0\rightarrow I\rightarrow F\mapright{\pi} A\rightarrow 0$
\end{center} 
constructed in proposition \ref{extens} and take the associated 
$C(X)$-extension 
\begin{center}
$0\rightarrow \K\ot I\ot\O2\rightarrow 
D=(\K\ot F\ot 1_\O2 )+(\K\ot I\ot\O2 )\rightarrow \K\ot A\rightarrow 0$. 
\end{center} 
The $C(X)$-nuclear quotient map $\sigma =\sigma\circ id_{\K\ot A}$ 
from the separable nuclear continuous field $\K\ot A$ to the quotient 
$D/(\K\ot I\ot\O2 )\subset
M(\K\ot I\ot\O2 )/(\K\ot I\ot\O2 )$ 
admits a $C(X)$-linear completely positive 
lifting $\K\ot A\rightarrow D$ ($\subset \K\ot [\O2\ot C(X)]\ot\O2$) 
thanks to proposition \ref{chnuc}. This means that the class of 
$\sigma$ is invertible in $Ext(X;\K\ot A,I\ot\O2 )$ (see the second 
remark following theorem \ref{homoto}). 

But the group $Ext(X;\K\ot A,I\ot\O2 ){}^{-1}$ is $C(X)$-linear 
homotopy invariant (theorem \ref{homoto}), hence zero since the 
endomorphism $\phi_2(a)=s_1as_1^*+s_2as_2^*$ of $\O2$ is homotopic to 
the identity map $id_{\O2}$ (\cite[proposition 2.2]{cuntz}) and 
so $[\theta ]=2[\theta ]$ in $Ext(X;\K\ot A,I\ot\O2 ){}^{-1}$ 
for any invertible extension $\theta$ of $\K\ot A$ by $I\ot\O2 $. 
As a consequence, the $C(X)$-extension defined by $\sigma$ is 
stably trivial. 
Furthermore, the identity representation of $D\subset M(\K\ot I\ot\O2 )$ 
is a continuous field of faithful representations (proposition 
\ref{extens}) and the assertion {\it 2.} of proposition~\ref{result} 
implies that the quotient morphism 
$(id_\K\ot\pi\ot id_\O2)$ from $D$ to $\K\ot A$ 
admits a cross section $\alpha$ which is a morphism of $C(X)$-algebras.

\medskip 
This monomorphism $\alpha$ is going to enable us to conclude by standard 
arguments, using theorem~\ref{kirchb} and the result of Elliott 
(\cite{elliott}) that the \cst-algebra $\O2$ is isomorphic to $\O2\ot\O2$. 

Choose a non-zero minimal projection $e_{11}$ in the 
\cst-algebra $\K$ of compact operators that we embed in $\O2$ and 
let $\varphi$ be a state on $\O2$ such that $\varphi (e_{11})=1$. 
If we take a unital completely positive projection $\psi$ of $\O2$ 
onto the nuclear \cst-subalgebra $A\subset\O2$ (theorem \ref{kirchb}), 
the composed map 
$$E=\left(\varphi\ot id_A\right)\circ 
\left( id_\O2 \ot [\pi\circ (\psi\ot id_{C(X)})]\ot\varphi\right)$$ 
is a unital $C(X)$-linear completely positive map 
from $\O2\ot [\O2\ot C(X)]\ot\O2$ onto $A$. 
Take also an isometry $u\in\O2\ot C(X)$ such that 
$\alpha (e_{11}\ot 1_A)=uu^*$ (corollary~\ref{unite}) and 
define the unital $C(X)$-algebra morphism 
\begin{center}
$\beta :A\longrightarrow \O2\ot [\O2\ot C(X)]\ot\O2\simeq\O2\ot C(X)$ 
\end{center}
by the formula $\beta (a)=u^*\alpha (e_{11}\ot a)u$. 
If $\widetilde{E} : \O2\ot C(X)\rightarrow A$ is the 
completely positive unital map $d\mapsto E(udu^*)$, 
one gets for all $a\in A$ the equality 
\begin{center} 
$(\widetilde{E}\circ\beta )(a)=(E\circ\alpha )(e_{11}\ot a)=a$
\end{center} 
\medskip
\noindent $\bullet$ 
Let us now come back to the general case of a compact space $X$. 

Define the continuous field $B=A\ot C([0,1])$ over the perfect compact 
space $Y={X\times [0,1]}$. According to the previous 
discussion, there exist a unital completely positive map 
$\widetilde{E} :\O2\ot C(Y)\rightarrow B$ and a $C(X)\ot C([0,1])$-linear 
monomorphism $\widetilde{\alpha} :B\rightarrow \O2\ot C(Y)$ such that 
$\widetilde{E}\circ\widetilde{\alpha} =id_B$. 
If $ev_1: C([0,1])\rightarrow\C$ is the evaluation map at 
$x=1\in [0,1]$, define the two maps 
$E:\O2\ot C(X)\rightarrow A$ and $\alpha :A\rightarrow \O2\ot C(X)$ by 
\begin{center}
$E(d)=(id_A\ot ev_1)\circ\widetilde{E} (d\ot 1_{C([0,1])} )$ and 
$\alpha (a)=(id_{\O2\ot C(X)}\ot ev_1)\circ\widetilde{\alpha} 
(a\ot 1_{C([0,1])})$. 
\end{center} 
Then $E$ is a unital $C(X)$-linear completely 
positive map, $\alpha$ is a unital $C(X)$-linear monomorphism and one has 
the identity $E\circ\alpha =id_A$. $\square$ 

\medskip\noindent {\bf Remark:} 
Assume that $X$ is a locally compact metrizable space 
and that the $C_0(X)$-algebra $A$ is a nuclear continuous field 
of \cst-algebras 
over $X$, where $C_0(X)$ denotes the algebra of continuous functions on 
$X$ vanishing at infinity. If $\widetilde{X}$ is the Alexandroff 
compactification of $X$, the unital $C(\widetilde{X} )$-algebra 
${\cal A}$ generated by $A$ and $C(\widetilde{X} )$ in the multiplier 
algebra $M[A\oplus C(\widetilde{X} )]$ is a separable unital continuous 
field of \cst-algebras over $\widetilde{X}$ 
(\cite[proposition 3.2]{bla1}). By theorem \ref{concl}, there exists 
therefore a $C(\widetilde{X} )$-linear monomorphism 
$\alpha :{\cal A}\hookrightarrow \O2\ot C(\widetilde{X} )$ and the 
$C_0(X)$-algebra $A$ is isomorphic to the $C_0(X)$-subalgebra $\alpha (A)$ 
of $\O2\ot C_0(X)$. 

\section{Concluding remarks}\label{avoir}
{\bf \ref{avoir}.1} 
A $C(X)$-subalgebra of $\O2\ot C(X)$ 
is by construction exact and continuous. Conversely, 
if $A$ is a non-zero exact separable unital continuous field of 
\cst-algebras over a perfect metrizable compact space $X$, one has 
by proposition \ref{extens} a $C(X)$-extension 
\begin{center}
$0\rightarrow I\rightarrow F\rightarrow A\rightarrow 0$
\end{center}
where $F$ is a $C(X)$-subalgebra of $\O2\ot C(X)$. 
If the identity map $A\rightarrow A=F/I$ admits a $C(X)$-linear 
completely positive lifting $A\rightarrow F$, the same method as the one 
used in theorem \ref{concl} will imply that the exact continuous field $A$ 
is isomorphic to a $C(X)$-subalgebra of the trivial continuous field 
${\O2\ot C(X)}$. 

It is therefore interesting to know whether this map admits 
a $C(X)$-linear completely positive lifting in the not discrete case. 

\bigskip\noindent {\bf \ref{avoir}.2} 
Let us have a look at one of the technical problems involved, 
the Hahn-Banach type extension property in the continuous 
field framework for finite type $C(X)$-submodules. 

\medskip 
Let $A$ be a separable unital continuous field of \cst-algebras over 
a compact metrizable space $X$ and let 
$D$ be a finitely generated $C(X)$-submodule which is an operator system. 
Assume that $\phi :D\rightarrow C(X)$ is a $C(X)$-linear unital 
completely positive map. Then for $x\in X$, there exists, thanks to 
\cite[proposition 3.13]{bla1}, a continuous field of states $\Phi^x$ 
on $A$, i.e. a $C(X)$-linear unital positive map from $A$ to $C(X)$, 
such that for all $d\in D$, $$\Phi^x (d)(x)=\phi (d)(x).$$ 
As a consequence, given $\varepsilon >0$ and a finite subset ${\cal F}$ 
of $D$, one can build by continuity and 
compactness a continuous field of states $\Phi$ on $A$ such that 
\begin{center}
$\max\{ \|\Phi (d)-\phi (d)\| ,d\in {\cal F}\} <\varepsilon$. 
\end{center}

But one cannot find in general any continuous field of states on $A$ 
whose restriction to $D$ is $\phi$. Indeed, consider the 
$C(\widetilde{\N} )$-algebra $A=\C^2\ot C(\widetilde{\N} )$ where 
$\widetilde{\N} =\N\cup\{\infty\}$ is the Alexandroff compactification 
of the space $\N$ of positive integers. 
Define the positive element $a\in C_\infty (\widetilde{\N} )A\subset A$ 
by the formulas 
\begin{center}
$a_n=\left\{ \begin{array}{ll}
(\frac{1}{n+1} ,0)&\mbox{if $n$ even}\\
(0,\frac{1}{n+1} )&\mbox{if $n$ odd}
\end{array}\right.$
\end{center}
and let $\phi$ be the $C(\widetilde{\N} )$-linear unital completely 
positive map with values in $C(\widetilde{\N} )$ defined on the 
$C(\widetilde{\N} )$-submodule generated by the two 
$C(\widetilde{\N} )$-linearly independent elements $1_A$ and $a$ 
through the formula 
\begin{center} 
${\phi (a)(n)=}\frac{1}{n+1}$ if $n<\infty$ and $\phi (a)(\infty )=0$. 
\end{center}
Suppose that the continuous field of states $\Phi$ is a 
$C(\widetilde{\N} )$-linear extension of $\phi$ to $A$. 
Then as A. Bauval already noticed it, one has the contradiction 
\begin{center}
$\begin{array}{rl}
1=\Phi (1_A)(\infty )&=\Phi\left( (1,0)\ot 1\right) (\infty ) + 
\Phi\left( (0,1)\ot 1\right) (\infty )\\
&=\lim_{n\rightarrow\infty} \Phi\left( (1,0)\ot 1\right) (2n+1) +
\lim_{n\rightarrow\infty} \Phi\left( (0,1)\ot 1\right) (2n)\\
&=0+0=0 .
\end{array}$
\end{center} 

\bigskip
\setcounter{section}{1}
\renewcommand{\thesection}{\Alph{section}}

\begin{center}
{\Huge {\bf Appendix}} {\Large {\it by Eberhard Kirchberg}}

{{\it (Humboldt Universit\"at zu Berlin)}}
\end{center}

In this appendix, we solve in proposition \ref{lift} the lifting question 
raised in paragraph 4.1 through a continuous generalisation of joint 
work of E.G. Effros and U. Haagerup on lifting problems for \cst-algebras 
(\cite{effha}, see also \cite{wa}). This result enables us to state 
the following characterisation of separable exact continuous fields 
of \cst-algebras: 
\begin{vth}\label{exact} 
Let $X$ be a compact metrizable space and $A$ be a (unital) 
separable continuous field of \cst-algebras over $X$. 

Then the \cst-algebra $A$ is exact if and only if there exists 
a (unital) monomorphism of $C(X)$-algebras $A\hookrightarrow\O2\ot C(X)$. 
\end{vth}

\medskip
Let us start with a technical $C(X)$-linear version of Auerbach's 
theorem (\cite[proposition 1.c.3]{lind}) for 
a continuous field of \cst-algebras $A$ over $X$ 
which gives us local bases over $C(X)$ with continuous 
coordinate maps for particular free $C(X)$-submodules of finite type 
in $A$. 

Define a $C(X)$--operator system 
in $A$ to be a $C(X)$-submodule which is an operator system. 
\begin{lemme}\label{auer} (\cite[lemma 2.4]{effha}) 
Let $A$ be a separable unital continuous field of \cst-algebras 
over a compact metrizable space $X$, $E\subset A$ be a a 
$C(X)$--operator system and assume that there exists an integer 
$n\in\N^*$ such that for all $x\in X$, the dimension $\dim E_x$ of the 
operator system $E_x\subset A_x$ equals $n$. Then the following holds. 

Given any point $x\in X$, there exist 
an open neighbourhood $\U$ of $x$ in $X$, 
self-adjoint $C(X)$-linear contractions $\varphi_i:A\rightarrow C(X)$ 
and self-adjoint elements $f_i\in E$ with $\| f_i\|\leq 2$ for 
$1\leq i\leq n$ such that 
\begin{center}
$\forall a\in C_0(\U)E$, $a=\sum_i \varphi_i(a)f_i$. 
\end{center}

Furthermore, there exists a continuous field of states 
$\Psi :A\rightarrow C(X)$ such that the restriction of the map 
$2n\Psi -id_A$ to the operator system $E$ is completely positive. 
\end{lemme}{\mbox{{\bf Proof :}}} 
Let us fix a point $x\in X$. Then there exist, thanks to 
Auerbach's theorem, a normal basis 
$\{ r_1,\ldots ,r_n\}$ of the 
fibre $E_x$ where each $r_i$ is self adjoint and norm one 
hermitian functionals $\phi_j:A_x\rightarrow\C$, $1\leq j\leq n$, 
such that $\phi_j(r_i)=\delta_{i,j}$. 

Consider the polar decomposition 
$\phi_j =\phi_j^+-\phi_j^-$ where $\phi_j^+$ and $\phi_j^-$ are 
positive functionals such that 
$1=\|\phi_j\| =\|\phi_j^+\| +\|\phi_j^-\|$. 
By \cite{bla1} lemme 3.12, there exist $C(X)$-linear positive maps 
$\varphi_j^+$ and $\varphi_j^-:A\rightarrow C(X)$ which extend the 
functionals $\phi_j^+$ and $\phi_j^-$ on the fibre $A_x$ 
to the $C(X)$-algebra $A$ with the property that 
$\varphi_j^+(1)=\|\phi_j^+\|$ and $\varphi_j^-(1)=\|\phi_j^-\|$. 
Take also $n$ norm $1$ self-adjoint elements $e_i\in E$ 
satisfying the equality $(e_i )_x=r_i$ and define the matrix 
$T=\left[\varphi_j(e_i )\right]_{i,j}\in M_n(\R )\ot C(X)$. 

One has by construction $T_x=1_{M_n(\R )}$; the set $\U_1\subset X$ 
of points $y\in X$ for which the spectrum of $T_y\in M_n(\R)$ is 
included in the open set $\{ z\in\C ,|z|>1/2\}$ is therefore an 
open neighbourhood of $x$ in $X$ (\cite[proposition 2.4 {\it b)}]{bla1}). 
If $\eta$ is a continuous function on $X$ with values in $[0,1]$ 
which is $0$ outside $\U_1$ and $1$ on an open neighbourhood $\U$ 
of the point $x\in X$, the self-adjoint elements 
$f_1,\ldots ,f_n$ of norm less than $2$ are then well defined 
in $C_0(\U_1)E$ by the formula 
\begin{center}
$T\left( \begin{array}{c}f_1\\ \vdots\\ f_n\end{array}\right) 
=\eta\left( \begin{array}{c}e_1\\ \vdots\\ e_n
\end{array}\right)$ 
\end{center}
and satisfy the relation $\varphi_j(f_i)(y)=\delta_{i,j}$ 
for each $y\in\U$ since the matrix $T_y$ is invertible, whence 
the desired equality for every $a\in C_0(\U )E$. 

\bigskip 
Keeping the same fixed point $x$, define now the continuous 
field of states $\Phi =\frac{1}{n}\sum_i (\varphi_i^++\varphi_i^-)$. 
Then one gets for all $a\in C_0(\U )E$ the equality: 
\begin{center}
$(2n\Phi -id_A)(a)={\displaystyle \sum_{1\leq i\leq n}}
\left[ \varphi_i^+(a)(2-f_i)+\varphi_i^-(a)(2+f_i)\right]$. 
\end{center}
The restriction of the map $(2n\Phi -id_A)$ to $C_0(\U )E$ is 
therefore completely positive and an appropriate partition of the 
unit $1_{C(X)}$ enables us to conclude. $\square$ 

\bigskip
Noticing that a $C(X)$-linear map $\sigma :A\rightarrow B$ between 
$C(X)$-algebras is completely positive if and only if each induced map 
$\sigma_x:A_x\rightarrow B_x$ is completely positive (see for instance 
\cite[proposition 2.9]{bla1}), 
the lemma \ref{auer} allows us to state a continuous version of 
theorem 2.5 of \cite{effha}. 
Replacing then the continuous field $A$ 
by $A\oplus M_{2^\infty}(\C )\ot C(X)$ 
(where $M_{2^\infty}(\C )=
\overline{\lim_{n\rightarrow\infty}M_{2^n}(\C )}$) 
and working with finitely generated $C(X)$--operator systems 
$E_k\subset A\oplus\bigcup_n M_{2^n}(\C )\ot C(X)$ for which the 
function $x\mapsto\dim (E_k)_x$ is continuous, one derives the 
following desired $C(X)$-linear completely positive lifting result. 

\begin{prop}\label{lift} (\cite[theorem 3.4]{effha}) 
Suppose that $A$ and $B$ are two unital separable exact 
continuous fields of \cst-algebras over a compact space $X$ with 
$A=B/J$ for some nuclear ideal $J$ in $B$. 

Then there exists a $C(X)$-linear unital completely positive 
lifting $A\rightarrow B$ of $id_A$.
\end{prop}{\mbox{{\bf Proof :}}} 
Let us define the two continuous fields 
${\cal A} =A\oplus M_{2^\infty}(\C )\ot C(X)$ and 
${\cal B} =B\oplus M_{2^\infty}(\C )\ot C(X)$. 
It is clearly enough to find a $C(X)$-linear unital completely positive 
cross section $\theta$ of the quotient morphism 
${\cal B}\rightarrow {\cal A}$ (by \cite[theorem 3.3]{bla1}). 

\medskip
Consider a dense sequence $\{a_k\}$ in the self-adjoint part of 
${\cal A}$ where each $a_k$ belongs to the dense subalgebra 
$A\oplus\bigcup_n M_{2^n}(\C )\ot C(X)$ of $\cal A$ and $a_1=1$. 
Let us show that we may assume inductively that $C(X)$--operator system 
$E_n$ generated by the $a_k$, $1\leq k\leq n$, satisfies the equality 
$\dim (E_n)_x=n$ for every $n\in \N^*$ and every $x\in X$. 
The inductive step is the following. Given $n\geq 2$, 
there exists by construction an integer $l$ 
such that $E_{n}\subset A\oplus M_{2^l}(\C )\ot C(X)$. Set 
$a_{n}'=a_{n}+2^{-n-1} d_l\ot 1_{C(X)}$ where 
\begin{center}
$d_l=1_{M_{2^l}(\C )}\ot 
\left( \begin{array}{cc}1&0\\ 0&0\end{array} \right) 
\in M_{2^{l+1}}(\C )\subset M_{2^\infty}(\C )$.
\end{center}
Then the $C(X)$-module 
$E_n'=E_{n-1}+C(X)a_n'$ verifies for each $x\in X$ the equality 
$\dim (E_n')_x =\dim (E_{n-1})_x +1$. 

Using proposition \ref{chnuc}, one can now finish the proof 
by the same method as the one developed by E.G. Effros and 
U. Haagerup in \cite{effha}.3 (see also \cite[theorem 6.10]{wa}). 
$\square$

\bigskip
\small Inst. Math. Luminy, CNRS--Luminy case 930, 
F--13288 Marseille CEDEX 9\\ 
\indent\indent E-mail: E.Blanchard@iml.univ-mrs.fr

\end{document}